\documentclass[12pt,a4paper]{article}
\usepackage{a4}
\usepackage{amssymb}
\usepackage{amsmath}
\usepackage{latexsym}
\usepackage{epsfig}
\usepackage{graphics}
\usepackage{color}

\newcommand{\A}{\ensuremath{\mathbb{A}}}
\newcommand{\Q}{\ensuremath{\mathbb{Q}}}
\newcommand{\Z}{\ensuremath{\mathbb{Z}}}

\begin{document}

\begin{center}
\bf MONOMIALIZATION OF MORPHISMS AND p-ADIC QUANTIFIER ELIMINATION  \\
\end{center}
\begin{center}
Jan Denef   \\
Department of Mathematics, K.U.Leuven, \\
Celestijnenlaan 200B, B-3001 Leuven, Belgium \\
\end{center}
We give a short proof of Macintyre's Theorem on Quantifier Elimination for the $p$-adic numbers, using a version of monomialization that follows directly from the Weak Toroidalization Theorem of Abramovich an Karu \cite{Abramovich1}(extended to non-closed fields \cite{Abramovich2}). \\ \\
\bf 1. Weak Monomialization of Morphisms \rm \\
\\
\bf Definition 1. \rm Let $f: X \rightarrow Y$ be a dominant morphism of nonsingular schemes of finite type over a field $K$ of characteristic zero, and let $D$ be a closed subset of $X$. The morphism $f$ is called \it monomial with respect to $D$ at \rm a $K$-rational point $a$ of $X$, if there exist an \'{e}tale morphism $U \rightarrow X$, a $K$-rational point $a'$ on $U$ mapping to $a$, uniformizing parameters $x_1, \dots, x_n$ on $U$, uniformizing parameters $y_1,\dots,y_m$ on a Zariski open neighborhood of $f(a)$ in $Y$, and a matrix $(a_{i,j})_{i,j}$ of non-negative integers (which necessarily has rank $m$) such that
$$ y_1 = x_1^{a_{11}} \cdots x_n^{a_{1n}},$$
(*) \qquad\qquad\qquad\qquad\qquad\qquad\qquad\qquad \vdots
$$ y_m = x_1^{a_{m1}} \cdots x_n^{a_{mn}},$$
and such that the preimage of $D$ in $U$ is contained in the locus of $\prod_i x_i$.
\\The morphism $f$ is called \it monomial with respect to \rm $D$ if, for each field extension $L$ of $K$, the morphism $f\otimes_K L : X \otimes_K L\, \rightarrow \, Y\otimes_K L\;$ is monomial with respect to the preimage of $D$ at each $L$-rational point of $X \otimes_K L$. The morphism $f$ is called \it monomial \rm if it is monomial with respect to the empty set.

With \emph{uniformizing parameters on} a scheme $Z$ over $K$, we mean a sequence of regular functions on $Z$ that induces an \'{e}tale morphism from $Z$ to affine space over $K$. When $A$ is a local ring containing its residue field, we will call elements $a_1, \dots, a_r \in A$ \emph{uniformizing parameters for} $A$, if these elements minus their residue are regular parameters for $A$.
\\ \\
\it Remark. \rm The definition of a monomial morphism as stated above is slightly different from Definition 1.2 in Cutkosky's paper \cite{Cutkosky1}, but is what we need for the present paper.
\\ \\
Next we state a version of monomialization that follows directly from the Weak Toroidalization Theorem of Abramovich an Karu \cite{Abramovich1}(extended to non-closed fields \cite{Abramovich2}). This version is also related to (but not implied by) Cutkosky's Theorem  on Local Monomialization of Morphisms (Theorem 1.3 in \cite{Cutkosky1}).
\\
\\
\bf Theorem 1 (Weak Monomialization). \it Let $f: X \rightarrow Y$ be a dominant morphism of integral separated schemes of finite type over a field $K$ of characteristic zero. Let $S$ be a proper closed subset of $X$. Then there exist nonsingular integral separated schemes $X_1$ and $Y_1$ of finite type over $K$, and birational proper morphisms $\alpha: X_1 \rightarrow X$, $\beta: Y_1 \rightarrow Y$, and a morphism $g: X_1 \rightarrow Y_1$ which is monomial with respect to $\alpha^{-1}(S)$, such that the following diagram commutes\\
$$X_1 \xrightarrow{\alpha }X $$
$$g \downarrow \quad \;\; \; \; \downarrow f$$
$$Y_1 \xrightarrow[\beta ]{} Y $$
and such that $\alpha^{-1}(S)$ is a strict normal crossings divisor on $X_1$.
\rm \\ \\
It is conjectured that we can take $\alpha$ and $\beta$ to be compositions of blow-ups of non-singular subvarieties. This conjecture is a weakening of the Conjecture on (Strong) Toroidalization of Morphisms \cite{Abramovich1, Cutkosky3}.
\\
\\
\it Remark. \rm To see that the Weak Toroidalization Theorem \cite{Abramovich2} implies Theorem 1, note that a dominant morphism $f: X \rightarrow Y$ of nonsingular varieties over $K$, which is toroidal with respect to toroidal embeddings $U_X \subset X, \; U_Y \subset Y$ (as defined in \cite{Abramovich2}), is also monomial with respect to $D_X = X \setminus U_X$. Indeed this follows from Proposition 1 in section 3, and the remark below it.
\\
\\
\bf 2. Quantifier elimination for $\Q_p$ \rm \\
\\
Let $p$ be a prime number. We denote the field of $p$-adic numbers by $\Q_p$.

Let $X$ be a scheme of finite type over $\Q_p$. A subset $A$ of $X(\Q_p)$ is called \it semi-algebraic \rm if $X$ can be covered by affine open subschemes $Y$ such that $A \cap Y(\Q_p)$ is a finite boolean combination of subsets determined by conditions of the form $f(y)= 0$, or of the form $\exists z \in \Q_p: f(y)=z^m$, with $f$ any regular function on $Y$ and $m$ any positive integer.

It is well known that any subset of $\Q_p^n$ of the form $C = \{y \in \Q_p^n \,|\, {\rm ord}_pf(y)\leq {\rm ord}_pg(y) \}$, with $f,g$ polynomials over $\Q_p$, is semi-algebraic. Indeed, for $p\neq 2$, $y \in C$  if and only if $(f(y))^2 + p(g(y))^2$ is a square in $\Q_p$, and for $p = 2$, $y \in C$ if and only if $(f(y))^2 + 8(g(y))^2$ is a square in $\Q_p$.

If $U$ is an open subscheme of $X$, then any semi-algebraic subset $B$ of $U(\Q_p)$ is a semi-algebraic subset of $X(\Q_p)$. Moreover, if $X$ is affine, then $B$ is expressible with regular functions on $X$, meaning that $B$ is expressible by a finite boolean combination of conditions of the form $f(y)= 0$, or $\exists z \in \Q_p: f(y)=z^m$, with $f$ any regular function on $X$. To verify this, note that this is evident when $X$ is affine,  $U$ is the complement of the locus of a regular function on $X$, and $B$ is expressible with regular functions $U$.

For any semi-algebraic subset $A$ of $X(\Q_p)$, there exist an affine reduced scheme $Z$ of finite type over $\Q_p$ and a morphism $h: Z \rightarrow X$, such that $A = h(Z(\Q_p))$. Indeed,  this follows easily from the observation that the complement of a coset of the subgroup of $m$-th powers in $\Q_p^\times$, is a finite union of such cosets.

We are now ready to state Macintyre's Theorem \cite{Macintyre1} on Quantifier Elimination for $\Q_p$, slightly reformulated. \\
\\
\bf Theorem 2 (Macintyre \cite{Macintyre1}). \it Let $f:X \rightarrow Y$ be a morphism of schemes of finite type over $\Q_p$. Let $A$ be a semi-algebraic subset of $X(\Q_p)$. Then $f(A)$ is a semi-algebraic subset of $Y(\Q_p)$.\rm \\
\\
Macintyre's proof is based on methods from mathematical logic and work of Ax and Kochen  \cite{Ax1} and Er\v sov \cite{Ershov1}. An elementary but rather convoluted proof has been given in \cite{Denef2}, based on $p$-adic cell decomposition which goes back to P. Cohen \cite{Cohen1}. For some applications to number theory we refer to \cite{Denef1}. We will now derive Macintyre's Theorem from Theorem 1 (weak monomialization) in a straightforward way. \\
\\
\it Proof of Theorem 2. \rm Clearly it suffices to prove the theorem when $X$, $Y$ are affine and integral, and $A=X(\Q_p)$. Let $\widetilde{X}, \widetilde{Y}$ be projective closures of $X, Y$. Replacing $Y$ by $\widetilde{Y}$, $X$ by the closure of the graph of $f$ in $\widetilde{X} \times \widetilde{Y}$, and $f$ by the projection morphism to $\widetilde{Y}$, we see that it suffices to prove the theorem when $X$, $Y$ are proper and integral, and $A = X(\Q_p) \setminus S(\Q_p)$, where $S$ is a proper closed subset of $X$. Moreover we may assume that $f$ is dominant.
Using Theorem 1, Lemma 1 below, and induction on dim $X$, we see that it is sufficient to prove the theorem when $X$ and $Y$ are proper, integral, and nonsingular, $f$ is dominant and \emph{monomial} with respect to a closed subset $D$ of $X$ of pure codimension 1, and $A = X(\Q_p) \setminus D(\Q_p)$. Then $X(\Q_p)$ is compact, and working locally, we reduce easily to the case that $X=\A^n, Y=\A^m$, $f$ is dominant and given by monomials in the standard coordinates, and $A$ is a $p$-adic ball minus the coordinate hyperplanes. Indeed, each point $a \in X(\Q_p)$ has arbitrary small p-adic neighborhoods that are balls in the local coordinates making $f$ monomial at $a$, and, using induction on dim $X$, we can cut out the coordinate hyperplanes. The theorem follows now from Lemma 3. \\
\\
\bf Lemma 1. \it Consider the following commutative diagram of morphisms of integral schemes of finite type over $\Q_p$
$$X_1 \xrightarrow{\alpha }X $$
$$g \downarrow \quad \;\; \; \; \downarrow f$$
$$Y_1 \xrightarrow[\beta ]{} Y $$
with $f$ dominant, and $\alpha$ and $\beta$ birational. Let $X_0 \varsubsetneqq X, Y_0 \varsubsetneqq Y$ be reduced closed subschemes with $f^{-1}(Y_0)\subset X_0$, such that $\alpha$ and $\beta$ are isomorphisms above $X\setminus X_0, Y\setminus Y_0$. Let $A$ be a subset of $X(\Q_p)$. If $f(X_0(\Q_p)\cap A)$ and $g(X_1(\Q_p)\cap \alpha^{-1}(A))$ are semi-algebraic, then $f(A)$ is semi-algebraic. \rm \\
\\
\it Proof. \rm Note that
$$ f(A)=\beta ( g(X_1(\Q_p)\cap \alpha^{-1}(A)) ) \cup f(X_0(\Q_p)\cap A), $$
because $\alpha$ is an isomorphism above $X\setminus X_0$. Thus
$$ f(A)=[ \beta ( g(X_1(\Q_p)\cap \alpha^{-1}(A)) )\setminus Y_0(\Q_p)] \cup f(X_0(\Q_p)\cap A), $$
because $f^{-1}(Y_0)\subset X_0$. The lemma now follows from the fact that for any semi-algebraic subset $B$ of $Y_1(\Q_p)$ we have that $\beta(B) \setminus Y_0(\Q_p)$ is semi-algebraic, because $\beta^{-1}(y)$ is given by rational functions in the coordinates of $y \in Y(\Q_p) \setminus Y_0(\Q_p) $, with no poles on $Y \setminus Y_0$. \\
\\
\newpage
\noindent
\bf Lemma 2. \it Let $K$ be a Henselian valued field. Fix elements $y_1,\cdots,y_m$ in $K^\times$, and consider the system of equations
$$ y_1 = x_1^{a_{11}} \cdots x_n^{a_{1n}},$$
\rm(*) \qquad\qquad\qquad\qquad\qquad\qquad\qquad\qquad \vdots
$$ y_m = x_1^{a_{m1}} \cdots x_n^{a_{mn}},$$
\it in the unknowns $x_1,\cdots,x_n$, with $m \leq n$, and  $A := (a_{i,j})_{i,j}$ a matrix of non-negative integers with rank equal to $m$. Let $d$ be the absolute value of the determinant of some {\rm (}$m,m${\rm )}-minor of $A$ with rank $m$. Let $I$ be any proper ideal of the valuation ring of $K$. Then any solution of {\rm (*)} in the group $K^\times /\,1+dI$ which can be lifted to a solution in $K^\times /\, 1+d^2I$, can be lifted to a solution of {\rm(*)} in the group $K^\times$. \rm \\
\\
\it Proof. \rm Let $x_1,\cdots,x_n \in K^\times$ be a solution of (*) in $K^\times /\, 1+d^2I$. Thus there exist $u_1,\cdots,u_m \in 1+d^2I$ such that the right side of (*) equals $y_1 u_1^{-1},\cdots,y_m u_m^{-1}$. We have to find $\varepsilon_1,\cdots,\varepsilon_n \in 1+dI$ such that
$x_1 \varepsilon_1,\cdots,x_n \varepsilon_n$ is an exact solution of (*). This is equivalent with the system of equations
$$ u_1 = \varepsilon_1^{a_{11}} \cdots \varepsilon_n^{a_{1n}},$$
(**) \qquad\qquad\qquad\qquad\qquad\qquad\qquad\qquad \vdots
$$ u_m = \varepsilon_1^{a_{m1}} \cdots \varepsilon_n^{a_{mn}}.$$
We may assume that $d=\det A_0$, where $A_0$ is the minor of $A$ consisting of the first $m$ columns of $A$. We choose $\varepsilon_{m+1}= \cdots = \varepsilon_n = 1$. Let $u$ be the column vector with components $u_1,\cdots,u_m$, and let $\varepsilon$ be the column vector with components $\varepsilon_1,\cdots,\varepsilon_m$. We have to find $\varepsilon \in (1+dI)^m$, such that $u=\varepsilon ^{A_0}$. A necessary condition is that $u^{B_0}=\varepsilon^{B_0 A_0}$, where $B_0 = dA_0^{-1}$ has integral coefficients, hence $u^{B_0}=\varepsilon^d$. Since $u \in (1+d^2I)^m$, we can indeed choose $\varepsilon \in (1+dI)^m$ such that $u^{B_0}=\varepsilon^d$. This implies $u^d=\varepsilon^{dA_0}$, hence $u=\varepsilon ^{A_0}$, because the components of both sides belong to $1+dI$ and have equal $d$-th powers. This concludes the proof of the lemma. \\
\\
\bf Lemma 3. \it Let $f: \Q_p^n \rightarrow \Q_p^m: (x_1,\cdots,x_n) \mapsto (y_1,\cdots,y_m) $ be given by the formula {\rm (*)} in Lemma 2, with $ (a_{i,j})_{i,j}$ a matrix of non-negative integers having rank $m$. Let $A \subset \Q_p^n$ be a $p$-adic ball minus the coordinate hyperplanes. Then f(A) is semi-algebraic. \rm \\
\\
\it Proof. \rm By covering and scaling, we may suppose that $A = (\Z_p \setminus \{0\})^k \times (1+I)^{n-k}$, with $I$ a proper ideal in $\Z_p$. By Lemma 2, an element $y \in \Z_p^m$ belongs to $f(A)$ if and only if (*) is solvable in $\Q_p^\times / 1+d^2I$ with side requirements $x_1, \cdots , x_k \in \Z_p$, and $x_{k+1}, \cdots ,x_n \in 1+I$. Clearly, this solvability condition on $y$ only depends on ord$_p \, y_1$, $\cdots$, ord$_p \, y_m$, and on the angular components modulo $d^2I$ of $y_1, \cdots, y_m$. The angular component of a $p$-adic number $z$ is defined as $zp^{-{\rm ord}_p \; z}$. Since there are a priori only a finite number of possibilities for the angular component modulo $d^2I$, the condition on the angular components modulo $d^2I$ of $y_1, \cdots, y_m$ is semi-algebraic. Indeed, it is easy to verify that the subset of $\Q_p$ consisting of all $p$-adic numbers with a fixed angular component modulo $d^2I$, is semi-albebraic (see e.g. \cite {Denef2}). Finally, by Presburger elimination, the condition on ord$_p \, y_1$, $\cdots$, ord$_p \, y_m$ is equivalent with a boolean combination of linear equations, linear inequalities, and linear congruences modulo fixed moduli. Again it is easy to verify that these are semi-algebraic, see e.g. \cite {Denef2}. This concludes the proof of Lemma 3.
\\
\\
\it Remark. \rm In the same way one can easily prove Tarski's elimination of quantifiers for the field of real numbers.
The above proof uses the local compactness of $\Q_p$ in an essential way. By a different argument we also obtained alternative proofs \cite{Denef3} of the Ax-Kochen-Er\v sov transfer principle and all classical results in the model theory of henselian valued fields of characteristic zero (except those on subanalytic sets). These alternative proofs are based on the Weak Toroidalization Theorem \cite{Abramovich1, Abramovich2}, and are purely algebraic geometric. To treat the theory of real and $p$-adic subanalytic sets in the same way, we would need a version of weak monomialization for analytic morphisms of analytic varieties.
\\
\\
\bf 3. Toroidal morphisms \rm \\
\\
In this section we prove a property of toroidal morphisms that we need in order to see that Theorem 1 is a direct consequence of the Weak Toroidalization Theorem \cite{Abramovich1, Abramovich2}.
\\  \\
\bf Proposition 1. \it Let $f: X \rightarrow Y$ be a dominant morphism of nonsingular integral separated schemes of finite type over a field $K$ of characteristic zero. Assume that $f$ is toroidal with respect to toroidal embeddings $U_X \subset X, \; U_Y \subset Y$ (as defined in \cite{Abramovich2}). Then for each $K$-rational point $a$ of $X$, there exist uniformizing parameters $\tilde{x}_1, \dots, \tilde{x}_n$ for the henselization $\mathcal{O}_{X,a}^h$ of $\,\mathcal{O}_{X,a}\,$, uniformizing parameters $y_1,\dots,y_m$ for $\mathcal{O}_{Y,f(a)}$, and non-negative integers $e_{j,i}\,$, such that $f$ is given by $y_j = \prod_i \tilde{x}_i^{\,e_{j,i}}$, the ideal of $D_X := X \setminus U_X$ in $\mathcal{O}_{X,a}^h$ is generated by $\prod_i \tilde{x}_i\,$, and the ideal of $D_Y := Y \setminus U_Y$ in $\mathcal{O}_{Y,f(a)}$ is generated by $\prod_j y_j$.
\rm \\
\\
\it Proof. \rm
Let $\overline{K}$ be an algebraic closure of $K$, and $\bar{X}$, $\bar{Y}$ the base change to $\overline{K}$ of $X$, $Y$. There exist uniformizing parameters $\bar{x}_i\,$, $i = 1, \dots , n$, for the completion of $\mathcal{O}_{\bar{X},a\,}$, and uniformizing parameters $\bar{y}_j\,$,  $j = 1, \dots , m$, for the completion of $\mathcal{O}_{\bar{Y},f(a)}$, such that $\bar{y}_j = \prod_{i = 1}^n \bar{x}_i^{\,e_{j,i}}$, with the $e_{j,i}$ forming a matrix of maximal rank $m \leq n$. Moreover, because $D_X$ and $D_Y$ are strict normal crossings divisors over $K$, there are uniformizing parameters $x_i\,$, $i = 1, \dots , n$, for $\mathcal{O}_{X,a\,}$, and uniformizing parameters $y_j\,$, $j = 1, \dots , m$, for $\mathcal{O}_{Y,f(a)}$, such that $\bar{x}_i/x_i,\; \bar{y}_j/y_j$ are units, the ideal of $D_X$ in $\mathcal{O}_{X,a}$ is generated by $\prod_{i=1}^n x_i$, and the ideal of $D_Y$ in $\mathcal{O}_{Y,f(a)}$ is generated by $\prod_{j=1}^m y_j$. Hence there exist units $u_j\,$, $j = 1, \dots , m$, in $\mathcal{O}_{X,a\,}$ such that
$y_j = u_j\prod_{i = 1}^n x_i^{e_{j,i}}$. Changing $y_j$ we may assume that $u_j$ has residue 1.

For $j = m+1, \dots , n$, define $u_j := 1$, and $y_j = u_j\prod_{i = 1}^n x_i^{\,e_{j,i}}$, with $e_{j,i}$ any non-negative integers such that the matrices
$$[e_{j,i}]_{j,i = 1, \dots n} \;\;\textrm{ and }\;\; [(\partial \log y_j / \partial \log x_i)(a)]_{j,i = 1, \dots n}$$
have maximal rank $n$. Such non-negative integers exist because the first $m$ rows of the second matrix form a minor with maximal rank $m$, since the same is true for $x_i,y_j$ replaced by $\bar{x}_i,\bar{y}_j$, for $j=1, \dots , m$, and because the components of the other rows of the second matrix are the $e_{j,i}$, with $j>m$.

By Hensel's Lemma there exist units $\epsilon_i \in \mathcal{O}_{X,a}^h\,$, with residue 1, for $i = 1, \dots , n$, such that $u_j = \prod_{i = 1}^n \epsilon_i^{e_{j,i}}$, for  $j = 1, \dots , n$. Put $\tilde{x}_i = \epsilon_i x_i \in \mathcal{O}_{X,a}^h\,$, for  $i = 1, \dots , n$, then $y_j = \prod_{i = 1}^n \tilde{x}_i^{\,e_{j,i}}$, for  $j = 1, \dots , n$. It remains to prove that the $\tilde{x}_i$ are uniformizing parameters for $\mathcal{O}_{X,a}^h$. Thus we have to show that the jacobian $\partial \tilde{x} / \partial x$ has maximal rank $n$ at $a$. But this follows directly from the formula
$\partial \log y / \partial \log x = \partial \log y / \partial \log \tilde{x} \;\, \partial \log \tilde{x} / \partial \log x$,
because the left side has maximal rank $n$ at $a$, and
$\det\partial \log \tilde{x} / \partial \log x = \left(\prod_i\epsilon_i^{-1}\right) \det \partial \tilde{x} / \partial x$.
This terminates the proof of the proposition.
\\
\\
\it Remark. \rm From the above proof it follows that in the definition of toroidal morphism \cite{Abramovich2}, we can replace the completions by henselizations. Using this fact, it is easy to verify that if $f$ is toroidal, then also $f\otimes_K L$ is toroidal for any field extension $L$ of $K$.
\\
\\
\it Acknowledgements. \rm We thank Dan Abramovich, Steven Dale Cutkosky, and Kalle Karu for stimulating conversations and information.

\end{document}